\documentclass[12pt]{amsart}
\usepackage{amssymb}

\numberwithin{equation}{section}

\newcommand{\ndash}{\nobreakdash-\hspace{0pt}}
\newcommand{\Ndash}{\nobreakdash--}

\newcommand{\ii}{{\mathrm{i}}}
\newcommand{\dd}{{\mathrm{d}}}

\newcommand{\id}{\mathrm{id}}
\newcommand{\ev}{\mathrm{ev}}

\DeclareMathOperator{\gh}{gh}

\newcommand{\frX}{\mathfrak{X}}
\newcommand{\frE}{\mathfrak{E}}

\newtheorem{Thm}{Theorem}[section]
\newtheorem{Prop}[Thm]{Proposition}
\newtheorem{Lem}[Thm]{Lemma}

\theoremstyle{remark}
\newtheorem{Rem}[Thm]{Remark}
\newtheorem*{Ack}{Acknowledgment}

\theoremstyle{definition}

\newtheorem{Exa}[Thm]{Example}

\newcommand{\braket}[2]{\left\langle{\,{#1}\,,\,{#2}\,}\right\rangle}

\newcommand{\Lie}[2]{{\left[{\,{#1}\,,\,{#2}\,}\right]}}

\newcommand{\BV}[2]{\left({\,{#1}\,,\,{#2}\,}\right)}

\newcommand{\sBV}[2]{\pmb({\,{#1}\,\pmb,\,{#2}\,}\pmb)}

\newcommand{\bbR}{{\mathbb{R}}}

\newcommand{\bbZ}{{\mathbb{Z}}}

\newcommand{\de}{\partial}

\newcommand{\Fun}{\mathit{Fun}}

\newcommand{\calC}{\mathcal{C}}

\newcommand{\calU}{\mathcal{U}}

\newcommand{\bomega}{\boldsymbol{\omega}}
\newcommand{\bvartheta}{\boldsymbol{\vartheta}}
\newcommand{\sfeta}{\boldsymbol{\eta}}
\newcommand{\bxi}{{\boldsymbol{\xi}}}
\newcommand{\bOmega}{\boldsymbol{\Omega}}

\newcommand{\sfA}{{\mathsf{A}}}

\newcommand{\sfS}{{\mathsf{S}}}

\newcommand{\sfe}{{\mathsf{e}}}
\newcommand{\sfX}{{\mathsf{X}}}
\newcommand{\sfQ}{{\mathsf{Q}}}

\newcommand{\frg}{{\mathfrak{g}}}

\DeclareMathOperator{\LL}{L}

\begin{document} 

\title{On the AKSZ formulation of the Poisson sigma model}

\author[A.~S.~Cattaneo]{Alberto~S.~Cattaneo}
\address{Institut f\"ur Mathematik, Universit\"at Z\"urich--Irchel,  
Winterthurerstrasse 190, CH-8057 Z\"urich, Switzerland}  
\email{asc@math.unizh.ch}

\author[G.~Felder]{Giovanni Felder}
\address{D-MATH, ETH-Zentrum, CH-8092 Z\"urich, Switzerland}
\email{felder@math.ethz.ch}
\subjclass{81T45 (Primary) 53D55, 58D30, 81T70 (Secondary)}
\keywords{Deformation quantization, topological quantum field theory,
Batalin--Vilkovisky (BV) formalism, $QP$\ndash manifolds}
\thanks{A.~S.~C. acknowledges partial support of SNF Grant
No.~2100-055536.98/1}

\begin{abstract}
We review and extend the Alexandrov--Kontsevich--Schwarz--Zaboronsky
construction of solutions of the Batalin--Vil\-ko\-vi\-sky classical
master equation. In particular, we study the case of sigma models
on manifolds with boundary. We show that a special case of this
construction yields the Batalin--Vilkovisky action functional 
 of the Poisson
sigma model on a disk. As we have shown in a previous paper, the
perturbative quantization of this model is related to Kontsevich's
deformation quantization of Poisson manifolds and to his formality
theorem. We also discuss the action of diffeomorphisms of the
target manifolds.

\end{abstract}

\maketitle

\section{Introduction}
In this paper we continue our study 
\cite{CF} of the Poisson sigma model on surfaces with boundary.   
In particular, we clarify the relation of our construction with   
the AKSZ method \cite{AKSZ} (see also \cite{S2}, \cite{P}
for more detailed descriptions), and comment on the   
action of diffeomorphisms.   
   
In \cite{CF} we discussed the quantization of Poisson sigma 
models---\hspace{0pt}topological sigma   
models whose target space is a Poisson manifold---on a disk in the   
framework of perturbative path integrals, and derived in this
context the
Kontsevich formula \cite{K} for deformation quantization \cite{BFFLS}.
   
The starting point   
is the classical action functional $S$, defined on the space of   
bundle maps from the tangent bundle $T\Sigma$ of a two-dimensional   
oriented manifold $\Sigma$, possibly with boundary,   
 to the cotangent bundle $T^*M$ of a Poisson   
manifold $M$. If we denote such a bundle map by a pair $(X,\eta)$,   
where $X\colon\Sigma\to M$ is the base map, and $\eta$, 
the map between fibers,   
 is a section   
in $\Gamma(\Sigma,\mathrm{Hom}(T\Sigma,X^*(T^*M)))=\Omega^1(\Sigma,X^*(T^*M))$,   
then the expression for the classical action functional is \cite{Ik,SS}   
\[   
S(X,\eta)=\int_\Sigma\langle\eta,dX\rangle+\frac12
(\alpha\circ X)(\eta,\eta).
\]   
Here the pairing is between the tangent and cotangent bundle of $M$ and the   
Poisson tensor $\alpha$ is viewed as a bilinear form on the
cotangent bundle.
The boundary conditions used in \cite{CF}
are $j_\partial^*\eta=0$ where 
$j_\partial\colon\partial\Sigma\to\Sigma$ denotes   
the inclusion of the boundary. 
Among the critical points of $S$ a special role is played by the {\em trivial   
classical solutions}\/ for which $X$ is constant and $\eta=0$. They are in one-to-one   
correspondence with points of $M$. The Feynman perturbation expansion around   
these trivial solutions    
was developed in the case of a disk $\Sigma$ in \cite{CF}. The gauge   
symmetry of the action functional was taken into account by applying   
the Batalin--Vilkovisky (BV) method \cite{BV} (see \cite{S} for a more precise 
mathematical description).
The Feynman rules for the resulting BV action functional   
yield, after gauge fixing, the Kontsevich formula 
for the star product \cite{K}.   
   
The AKSZ method \cite{AKSZ} is a method to construct solutions of the   
BV master equation directly, without starting from   
a classical action with a set of symmetries, as is done in the 
BV   
method. The classical action is then recovered {\em a posteriori}\/   
by setting the fields of non-zero degree to zero.   
This method is applied in \cite{AKSZ} to some examples and it is shown that   
the BV  action of the Chern--Simons theory, the Witten    
A- and B-models are special cases of their construction.   
   
Here we adapt the AKSZ method to the case of manifolds with boundary, and   
show that the BV action of the Poisson sigma model
derived in \cite{CF}, with the same boundary conditions, can be obtained
in this framework. 
Closely related results have been obtained recently in \cite{P}.  

In the case of the Poisson sigma model,
the construction goes as follows: to a
an oriented two-dimensional manifold $\Sigma$, possibly with boundary, one
associates the supermanifold $\Pi T\Sigma$, the tangent bundle
with reversed parity of the fiber. The algebra of functions on
$\Pi T\Sigma$ is isomorphic to 
the (graded commutative) algebra of differential forms 
on $\Sigma$ with values in a graded commutative ground ring $\Lambda$, which
we take, following \cite{V}, to be the infinite-dimensional exterior
algebra $\varinjlim\wedge\mathbb{R}^k$.
The integration of differential forms and the de~Rham differential
are, in the language of supermanifolds, a measure $\mu$ and a self-commuting
vector field $D$, respectively, on $\Pi T\Sigma$.
Let $M$ be a Poisson manifold and $\Pi T^*M$ be the cotangent bundle of $M$
with reversed parity of the fiber. This supermanifold has a canonical
symplectic structure (with exact symplectic form) just as any cotangent
bundle. It is odd, because of the parity reversion, meaning that the
corresponding Poisson bracket is graded skew-symmetric only after shifting the
degree by one. The algebra of functions on $\Pi T^*M$ is
isomorphic to the algebra of multivector fields on $M$ 
with values in $\Lambda$.
In particular, the Poisson bivector field on $M$ can be identified with
a function on $\Pi T^*M$. Its Hamiltonian vector field is an odd
vector field $Q$. It commutes with itself owing to the Jacobi identity.
The space of maps $\Pi T\Sigma\to\Pi T^*M$ inherits then an
odd symplectic form and an odd Hamiltonian self-commuting vector field:
the symplectic form is obtained from the symplectic form on $\Pi T^*M$
upon integration over $\Pi T\Sigma$ with respect to the measure $\mu$.
The odd vector field is the sum of the commuting vector fields $\hat D$ and 
$\check Q$ obtained from $D$ and $Q$ by acting on maps 
$\Pi T\Sigma\to\Pi T^*M$ by the corresponding
infinitesimal diffeomorphisms of $\Pi T\Sigma$ on the left and
of $\Pi T^*M$ on the right. The Hamiltonian function of this vector
field is then the BV action functional, if $\Sigma$
has no boundary. It is a function on
the space of maps from $\Pi T\Sigma$ to $\Pi T^*M$, whose Poisson
bracket with itself vanishes, i.e., it solves the BV
classical master equation. If $\Sigma$ has a boundary,
suitable boundary conditions must be 
imposed. They are first-class constraints, and the correct action
functional is obtained after Hamiltonian reduction.

We describe this construction in detail in Sect.~\ref{s-3} after reviewing 
and extending the AKSZ construction in Sect.~\ref{s-2}.
In the first part of Sect.~\ref{s-2}, we present the original construction
of \cite{AKSZ}, see Theorem \ref{thm-Dinv}.
In the typical case,
the data are a closed manifold $\Sigma$ and  an odd or even
(depending on the dimension of $\Sigma$) 
symplectic supermanifold $Y$ with a Hamiltonian odd self-commuting vector
field $Q$. The construction
produces a solution of the classical BV master equation 
on the space of maps from $\Pi T\Sigma$ to $Y$. 
In the case where $\Sigma$ has a boundary, a similar result holds if
the symplectic form on $Y$ is exact, and it depends on the choice of 
a symplectic potential---a $1$-form whose differential is the symplectic
form. Moreover the Hamiltonian function of $Q$ must vanish on the zero set
of the symplectic potential. Under these hypotheses, we construct a 
BV action functional on a Hamiltonian reduction of
the space of maps from $\Pi T\Sigma$ to $Y$ obeying suitable boundary conditions,
see Theorem \ref{thm-dementor}.
In the case when $Y$ has global Darboux coordinates, 
similar results---but with a different treatment of
boundary conditions---were obtained in \cite{P}.

Along the way, we also describe some by-products, such
as a geometric explanation of the supersymmetry 
discovered in \cite{CF} and generalized in \cite{P}.
We also comment on the action of diffeomorphisms of the target manifold $M$
of the Poisson sigma model. We show that 
diffeomorphisms of  $M$ induce canonical transformations
of the space of maps from $\Pi T\Sigma$ to $\Pi T^*M$. This means that one
should expect
the partition function to be invariant under such diffeomorphisms. However,
when we introduce boundary observables in the 
path integrals, the variation of the
correlation functions under a canonical transformation is singular, and
the diffeomorphism invariance is spoiled after regularization, 
albeit in a controlled way. These ideas, closely related to ideas in
\cite{L}, are a motivation for the globalization results in the 
form described in \cite{CFT}.

\begin{Ack}
We acknowledge an invitation to the 
Conf\'erence Mosh\'e Flato 2000 and especially thank Daniel
Sternheimer for his efforts in making the event successful and pleasant.
We also thank Carlo Rossi for stimulating discussions and useful comments.
\end{Ack}

\section{The AKSZ construction}\label{s-2}
The Batalin--Vilkovisky (BV) method in quantum field theory \cite{BV}
relies on extending the space of fields of a theory to a larger superspace
with a $QP$-structure (in the terminology introduced in \cite{S} and 
\cite{AKSZ}). The AKSZ method \cite{AKSZ} is a way of inducing a
$QP$-structure on a space of supermaps (given certain structures on
the domain and on the target supermanifolds). 
We begin by giving a short description of $Q$-, $P$- and $QP$-structures.

\subsection{Some structures on differentiable supermanifolds}
Let $N$ be a differentiable supermanifold. 
We will denote by $\Fun(M)$ its $\bbZ_2$\ndash graded algebra 
of smooth functions. Following
\cite{V}, we view $\Fun(M)$ as a $\Lambda$\ndash module,
where
\[
\Lambda:=\varinjlim\wedge^\bullet\bbR^k
\]
is the exterior algebra of $\bbR^\infty$.

\subsubsection{$Q$-structures}
A {\sf $Q$\ndash structure} on a supermanifold
is the choice of an odd 
self-commuting vector field. In homogeneous 
local coordinates $\{y_1,\dots,y_n\}$,
\[
Q=Q^i\de_i,
\]
with $\de_i=\de/\de y^i$.
``Odd" means that $\deg Q^i\equiv\epsilon_i+1\mod2$, where by $\epsilon_i$
we denote the degree of the coordinate $y^i$. ``Self-commuting" means
\[
\Lie QQ=2Q^2=2(Q^j\de_jQ^i)\de_i=0.
\]
In other words, a $Q$-structure is the choice of a differential
on the algebra of functions of $N$.

\begin{Exa}\label{exa-ptsigma}
Let $N=\Pi T\Sigma$, with $\Sigma$ an ordinary manifold. 
By this we denote the tangent bundle of $\Sigma$ with reversed parity on the 
fiber. By definition the algebra of functions on $\Pi T\Sigma$ is isomorphic
to the algebra of differential forms on $\Sigma$ tensor $\Lambda$.
Let us denote by $\phi$ this isomorphism.
Then the exterior derivative 
$\dd$ defines a $Q$-structure: $Df:=\phi^{-1}(\dd\phi(f))$,
$f\in\Fun(\Pi T\Sigma)$. Choosing local coordinates
$\{u^1,\dots,u^s\}$ on $\Sigma$ together with their odd counterparts
$\{\theta^1,\dots,\theta^s\}$, we can write
\[
D=\theta^\mu\frac{\de}{\de u^\mu}.
\]
Another example of differential is the contraction $\iota_v$ by a vector
field $v$. We define then the $Q$\ndash vector field
$K_v:=\phi^{-1}\iota_v\phi$ and
can write 
\[
K_v=v^\mu(u)\frac\de{\de \theta^\mu}
\]
in local coordinates. Finally, observe that the Lie bracket
$\Lie D{K_v}=:L_v$ is an even vector field corresponding
to the Lie derivative $\LL_v$ via the formula $L_v=\phi^{-1}\LL_v\phi$.
\end{Exa}

\subsubsection{$P$-structures}
A {\sf symplectic structure} is a closed nondegenerate $2$\ndash form on 
$N$. With respect to the $\bbZ_2$\ndash grading of $\Fun(M)$,
it can be even or odd. In local coordinates we write
\[
\omega=\frac12\dd y^i\omega_{ij}\dd y^j,
\]
with $\wedge$-symbols suppressed and the Koszul sign rule 
\begin{equation}\label{Kos}
\dd y^i\dd y^j=-(-1)^{\epsilon_i\epsilon_j}\dd y^j\dd y^i.
\end{equation}
An even symplectic structure satisfies
$\deg\omega_{ij}\equiv\epsilon_i+\epsilon_j$, while an odd symplectic
structure (shortly a {\sf $P$\ndash structure}) satisfies 
$\deg\omega_{ij}\equiv\epsilon_i+\epsilon_j+1$. 
Moreover,
\[
\omega_{ji}=\begin{cases}
(-1)^{(\epsilon_i+1)(\epsilon_j+1)}\,\omega_{ij} &\text{even case,}\\
(-1)^{1+\epsilon_i\epsilon_j}\,\omega_{ij} &\text{odd case.}
\end{cases}
\]
A symplectic structure associates to each function $f$ a vector field 
$X_f$ by the identity $\iota_{X_f}\omega=\dd f$. 
Observe that $X_f$ has the same (opposite) parity of $f$ if $\omega$ is
even (odd).
One defines the bracket
of two functions $f$ and $g$ by $\BV fg:=X_f(g)=\iota_{X_f}\iota_{X_g}\omega$. 
In the even case,
it has the properties of a Poisson bracket, while in the odd case it is
a Gerstenhaber bracket (but is referred to in the quantum field theory 
literature as a BV bracket).

\begin{Exa}\label{exa-ptm}
Let $N=\Pi T^* M$, with $M$ an ordinary manifold.
The algebra of functions on $N$ can be identified
with the algebra of multivector fields 
on $M$ tensor the infinite-dimensional exterior algebra $\Lambda$.
The algebra of multivector fields admits
a nondegenerate Gerstenhaber bracket: viz., the Schouten--Nijenhuis bracket
$\Lie\ \ $. This bracket is determined by the canonical odd symplectic
form $\omega$ that, using local coordinates $x^1,\dots x^m$ on $M$
and their odd counterparts $p_1,\dots,p_m$, can be written as
$\omega=\dd p_i\,\dd x^i$.
\end{Exa}

\subsubsection{$QP$-structures}
A vector field $Q$ and a symplectic structure are said to be {\sf compatible}
if $\LL_Q\omega=\dd\iota_Q\omega=0$. 
Here and in the following, $\LL_X=\dd\iota_X+\iota_X\dd$ denotes
the Lie derivative w.r.t.\ the (odd or even) vector field $X$.
If in particular $\iota_Q\omega=\dd S$
for some function $S$, then $Q$ is said to be {\sf Hamiltonian}.
A {\sf $QP$\ndash manifold} is
a supermanifold with compatible $Q$- and $P$\ndash structures; in case
$Q$ is Hamiltonian, its Hamiltonian function $S$ is even and satisfies
the so-called master equation $\BV SS=0$.

\begin{Exa}\label{exa-ptmPQ}
Again let $N=\Pi T^* M$. Let $\alpha$ be a multivector field on $M$
of even degree satisfying $\Lie\alpha\alpha=0$. Then $\alpha$
can be seen as the Hamiltonian function of an odd self-commuting vector
field $Q$ that acts on multivector fields by $\Lie\alpha\ $.
For example, if $\alpha$ is a Poisson bivector field,
we may write in local coordinates
\begin{subequations}\label{SQ}
\begin{align}
S_\alpha(x,p)&=\frac12\alpha^{ij}(x)p_ip_j,\label{SQ-S}\\
Q_\alpha(x,p)&=\alpha^{ij}(x)p_j\frac\de{\de x^i}+
\frac12\de_i\alpha^{jk}p_jp_k\frac\de{\de p_i}.\label{SQ-Q}
\end{align}
\end{subequations}
\end{Exa}

\subsubsection{Measures}
We recall that a {\sf measure} on an $(m,n)$\ndash supermanifold $N$ 
is a linear functional on the algebra of function $\Fun(M)$ that
kills all components of homogeneous degree in the odd
coordinates less than $n$. More precisely, a measure is a section
of the Berezinian bundle.
Our notation for a measure applied to a function $f$ is
$\int_N f\;\mu$, and we will often call $\mu$ the measure.
We say that a measure is {\sf nondegenerate} if
its composition with the 
product yields a nondegenerate bilinear form on $\Fun(M)$.
Finally, given a vector field $D$ on $N$, we say that the measure
is $D$\ndash{\sf invariant} if $\int_N Df\;\mu=0$, $\forall f\in\Fun(N)$.

\begin{Exa}\label{exa-mu}
Let $N=\Pi T\Sigma$ as in Example~\ref{exa-ptsigma}.
Its {\sf canonical  measure} $\mu$  is defined by
$\int_{\Pi T\Sigma}f\,\mu=\int_\Sigma\phi(f)$,
$f\in\Fun(\Pi T\Sigma)$, and is clearly nondegenerate.
In local coordinates, assuming the rules
of Berezinian integration ($\int\theta^\mu\,\dd\theta^\nu=\delta^{\mu\nu}$), 
we may write
\[
\mu=\dd\theta^s\dotsm\dd\theta^1\,\dd x^1\dotsm\dd x^s.
\]
This measure is compatible with the vector field $D$ defined
in Example~\ref{exa-ptsigma} in the following sense.
If $\Sigma$ has a boundary $\de\Sigma$, we denote by $\mu^\de$
the canonical measure on $\Pi T\de\Sigma$ and by $ i_\de$ the inclusion
$\Pi T\de\Sigma\hookrightarrow\Pi T\Sigma$ induced from
$\de\Sigma\hookrightarrow\Sigma$.\footnote{More precisely, let
$j_\de$ be the inclusion of $\de\Sigma$ into $\Sigma$. 
Then $i_\de^*\colon\Fun(\Pi T\Sigma)\to\Fun(\Pi T\de\Sigma)$ is
defined as $\phi_\de^{-1}(j_\de^*\otimes\id)\phi$,
where $\phi\colon\Fun(\Pi T\Sigma)\to\Omega^\bullet(\Sigma)\otimes\Lambda$
and $\phi_\de\colon\Fun(\Pi T\de\Sigma)\to\Omega^\bullet(\de\Sigma)\otimes\Lambda$
are the isomorphisms that define the functions on these supermanifolds.}
Then Stokes' theorem can be reformulated as follows:
\begin{equation}\label{intDf}
\int_{\Pi T\Sigma}Df\;\mu = \int_{\Pi T\de\Sigma} i_\de^*f\;\mu^\de.
\end{equation}
In particular, if $\de\Sigma=\emptyset$, the canonical measure
$\mu$ is $D$\ndash invariant.
Finally, for every vector field $v$ on $\Sigma$, the canonical
measure $\mu$ is $K_v$\ndash invariant, with $K_v$ defined in 
Example~\ref{exa-ptsigma}, since the top form component of $\iota_v\phi(f)$
vanishes for every $f\in\Fun(\Pi T\Sigma)$.
\end{Exa}

More generally, let $N$ and $L$ be two supermanifolds and let
$\mu$ be a measure on $N$. We may define a chain map
$\mu_*\colon\Omega^\bullet(N\times L)\to\Omega^\bullet(L)$, $\forall k$, by the rule
\[
(\mu_*\omega)(z)(\lambda_1,\dots,\lambda_k)=
\int_{y\in N}\omega(y,z)(\lambda_1,\dots,\lambda_k)\;\mu(y),
\]
with $z\in L$ and $\lambda_1,\dots,\lambda_k\in T_zL$.
Observe that, if $N$ is of type $(m,n)$ with $n$ odd (even),
then $\mu_*$ changes (preserves) the parity.
Moreover, if $\mu$ is $D$\ndash invariant, one obtains
$\mu_*\LL_{D_1}=0$, where $D_1$ is the lift to $N\times L$ 
of the vector field $D$ defined on $N$.

In the case $N=\Pi T\Sigma$ described in Example~\ref{exa-mu},
the generalization of \eqref{intDf} is instead
\begin{equation}\label{muLD}
\mu_*\LL_{D_1}= \mu^\de_*(\iota_\de\times\id)^*,
\end{equation}
where $\mu$ is the canonical measure, and $D$ is the vector
field defined in Example~\ref{exa-ptsigma}.

\subsubsection{$\bbZ$-grading}\label{ssec-zgrad}
In the application of the BV formalism to quantum field theories, one 
usually considers a $\bbZ$\ndash grading ({\sf ghost number}) instead
of just a $\bbZ_2$-grading. This is obtained as follows. Let $A$ be
a $\bbZ_2$\ndash graded supercommutative algebra (e.g.,
the algebra of functions on some supermanifold). Define
$\sfA:=A[\bbZ]^{\bbZ_2}$,
where the action of the nontrivial element $\epsilon\in\bbZ_2$ on 
$n\in\bbZ$ is given by $\epsilon n=(-1)^n n$, while its action on $A_+$ 
is trivial and $\epsilon a=-a$ for $a\in A_-$.
We can decompose $\sfA$ as $\bigoplus_{j\in\bbZ}\sfA_j$ with
\begin{align*}
\sfA_{2j} &= A_+[2j],\\
\sfA_{2j+1} &= A_-[2j+1].
\end{align*}
We say then that an element of $\sfA_j$ has ghost number $j$.
It follows that $\sfA$ is a $\bbZ$\ndash graded supercommutative algebra.

We finally have to define an inclusion of the original algebra 
$A$ into $\sfA$. Consider the case when $A$ is the algebra of functions
on a supermanifold of the form $\Pi E$, with $E$ a vector bundle over
an ordinary manifold:
$A=\Fun(E)\simeq\Gamma(\wedge^\bullet E^*)\otimes\Lambda$.
We define then the inclusion by sending the even coordinates
to $\sfA_0$ and the linear
odd coordinates (i.e., those corresponding to linear 
functions on the fiber of $E$)
to $\sfA_1$. (That is, we assign
ghost number zero to the former and one to the latter.)

More generally, if $A=\Omega^\bullet(N)$ is the algebra of differential
forms on a supermanifold $N=\Pi E$, the inclusion is defined by sending
even coordinates and their differentials to $\sfA_0$ and 
linear odd coordinates
and their differentials to $\sfA_1$.
Then $\bOmega^\bullet(N)$ is a bigraded algebra and we follow the Koszul
sign rule \eqref{Kos} for the products.
As a consequence of this, the chain map $\mu_*$ defined above
can be extended to a chain map 
$\bOmega^\bullet(N\times L)\to\bOmega^\bullet(L)$
that lowers the ghost number by $n$ 
if $N=\Pi E$ with $E$ a vector bundle of rank $n$.

Observe that the vector fields $D$ and $K_v$ defined in 
Example~\ref{exa-ptsigma} have ghost number $1$ and $-1$ respectively,
while the vector field $Q$ in \eqref{SQ-Q} has ghost number $1$ being
generated by the ghost-number-zero function $S$ in \eqref{SQ-S} using
the canonical symplectic form $\omega$ of Example~\ref{exa-ptm}
that has ghost number $1$.

In this framework, one usually speaks of a $Q$\ndash structure only
if the self-commuting vector field $Q$ has ghost number $1$ 
and of $P$\ndash structure only if the symplectic form
$\omega$ has ghost number $-1$ 
(i.e., $\gh Q^i=\epsilon_i+1$ 
and $\gh\omega_{ij}=-\epsilon_i-\epsilon_j-1$,
with $\epsilon_i$ the ghost number of the $i$th coordinate).
So the Hamiltonian function $S$ of $Q$, when existing,
has ghost number $0$ and is usually called the {BV\ndash action}.

\subsection{Spaces of maps}
The AKSZ construction is a way of defining a $QP$-structure on the
space $Y^X$ of smooth maps from the supermanifold $X$ to the 
supermanifold $Y$.

\subsubsection{$Q$-structure}\label{subs-Q}
Diffeomorphisms of $X$  and of $Y$ act on $Y^X$. As the former act
from the right and the latter from the left, they (super)commute.
At the infinitesimal level, this means that we can associate to
vector fields $D$ on $X$ and $Q$ on $Y$ the (super)commuting vector
fields $\Hat D$ and $\Check Q$ on $Y^X$. 
Recall that $T_f Y^X\simeq\Gamma(X,f^*TY)$, $f\in Y^X$. So
a vector field on $Y^X$ assigns to each $x\in X$ and each $f\in Y^X$
an element of $T_{f(x)}Y$. So we can write
\[
\Hat D(x,f)=\dd f(x)D(x),\qquad
\Check Q(x,f)=Q(f(x)).
\]
In particular, if
both $D$ and $Q$ are odd and self-commuting, any linear combination
of $\Hat D$ and $\Check Q$ defines a $Q$-structure on $Y^X$.

\subsubsection{$Y$-induced $QP$-structure}
Let $\mu$ be a nondegenerate measure on $X$ and $\omega$ a symplectic 
structure on $Y$. 
We assume $\omega$ to be even (odd) if $X$ is of type $(m,n)$
with $n$ odd (even). Consider the evaluation map
\[
\ev\colon\begin{array}[t]{ccc}
X\times Y^X &\to& Y,\\
(x,f) &\mapsto& f(x).
\end{array}
\]
This induces a chain map $\mu_*\ev^*\colon\bOmega^\bullet(Y)\to\bOmega^\bullet(Y^X)$
that lowers the ghost number by $n$ if $X$ is of type $(m,n)$.
Then $\bomega:=\mu_*\ev^*\omega$ defines a $P$\ndash structure on
$Y^X$. 

It is easy to verify that $\Check Q$ is compatible with $\bomega$.
In fact, $\iota_{\Check Q}\bomega=\mu_*\ev^*\iota_Q\omega$.
In particular, if $Q$ is Hamiltonian with Hamiltonian function $S$,
then so is $\Check Q$ with Hamiltonian function
\begin{equation}\label{checks}
\Check S=\mu_*\ev^*S.
\end{equation}
In fact, 
\[
\dd\Check S=\mu_*\ev^*\dd S=
\mu_*\ev^*\iota_Q\omega=\iota_{\Check Q}\bomega.
\]
In the following, we will denote by $\sBV\ \ $ the Gerstenhaber bracket
induced by $\bomega$ on $\Fun(Y^X)$.
\begin{Prop}
$\mu_*\ev^*\colon(\Fun(Y),\BV\ \ )\to(\Fun(Y^X),\sBV\ \ )$ is a 
Lie algebra homomorphism.
\end{Prop}
\begin{proof}
Let $F$ and $G$ be two functions on $Y$, and denote by $Q$
the Hamiltonian vector field of $F$. So $\BV FG=\iota_Q\dd G$.
As already observed,
$\Check Q$ is the Hamiltonian vector field of $\mu_*\ev^*F$.
As a consequence,
\[
\sBV{\mu_*\ev^*F}{\mu_*\ev^*G}=\iota_{\Check Q}\dd\mu_*\ev^*G=
\mu_*\ev^*\iota_Q\dd G=\mu_*\ev^*\BV FG.
\]
\end{proof}

\subsubsection{The AKSZ $QP$-structure}
We have constructed above a $P$\ndash structure $\bomega$
on $Y^X$ as well as a $Q$\ndash structure $\sfQ:=\Hat D+\Check Q$.
We have also checked that $(\Check Q,\bomega)$ is a $QP$\ndash structure.
What is left to check is whether $\Hat D$ is also compatible with
$\bomega$. We need the following
\begin{Lem}\label{lem-Dinv}
If $\mu$ is $D$\ndash invariant, then $\LL_{\Hat D}\mu_*\ev^*=0$.
\end{Lem}
\begin{proof}
Since the evaluation map is invariant under the simultaneous action 
of the inverse of a diffeomorphism on $X$ and of the induced
diffeomorphism on $Y^X$---i.e., $\ev(\phi^{-1}(x),f\circ\phi)=\ev(x,f)$)---, 
we obtain, at the infinitesimal level,
\[
\LL_{D_1}\ev^*=\LL_{\Hat D_2}\ev^*,
\]
where $D_1$ and $\Hat D_2$ are the lifts of $D$ and $\Hat D$ to
$X\times Y^X$. On the other hand, by definition of $\mu_*$,
we have $\iota_{\Hat D}\mu_*=\mu_*\iota_{\Hat D_2}$.
Combining these identities, we get
\[
\LL_{\Hat D}\mu_*\ev^*=\mu_*\LL_{\Hat D_2}\ev^*=
\mu_*\LL_{D_1}\ev^*=0.
\]
\end{proof}
{}From this we obtain the
\begin{Thm}\label{thm-Dinv}
If $D$ and $Q$ are $Q$\ndash structures on $X$ and $Y$ respectively,
$\mu$ is a nondegenerate $D$\ndash invariant measure on $X$ and
$\omega$ is a symplectic form on $Y$ (of parity opposite to the
odd dimension of $X$), then $(\Hat D+\Check Q,\bomega)$ is
a $QP$\ndash structure on $Y^X$.
If moreover $\omega=\dd\vartheta$, then $\Hat D$ is Hamiltonian with
Hamiltonian function
\[
\Hat S=-\iota_{\Hat D}\bvartheta,
\]
with $\bvartheta=\mu_*\ev^*\vartheta$.
If in addition $Q$ admits a Hamiltonian function $S$, then
$\sBV{\Hat S}{\Check S}=0$, with $\Check S$ defined in \eqref{checks}.
In particular $\sfS=\Hat S+\Check S$ satisfies the master equation
$\sBV\sfS\sfS=0$.
\end{Thm}

\subsubsection{The case $X=\Pi T\Sigma$}\label{ssec-ptsigma}
As observed in Example~\ref{exa-mu}, the canonical measure $\mu$
on $\Pi T\Sigma$ is $D$\ndash invariant when $D$ is the vector
field defined in Example~\ref{exa-ptsigma} and $\de\Sigma=\emptyset$.
This is the typical case in which the AKSZ formalism is applied.
We want however  to consider also the case when $\Sigma$ has a boundary.
In this case we cannot apply Lemma~\ref{lem-Dinv}. 
We have however the following generalization.
\begin{Lem}\label{lem-ldhat}
Let $\mu^\de$ be the canonical measure on $\Pi T\de\Sigma$,
and $ i_\de$ the inclusion $\de\Sigma\hookrightarrow\Sigma$.
Define $\ev_\de^*= (i_\de\times\id)^*\ev^*\colon\bOmega^\bullet(Y)
\to\bOmega^\bullet(\Pi T\de\Sigma\times Y^{\Pi T\Sigma})$.
Then 
\[
\LL_{\Hat D}\mu_*\ev^*=\mu^\de_*\ev_\de^*.
\]
\end{Lem}
\begin{proof}
Proceed as in the proof of Lemma~\ref{lem-Dinv} till the last equation
which is now replaced by
\[
\LL_{\Hat D}\mu_*\ev^*=\mu_*\LL_{\Hat D_2}\ev^*=
\mu_*\LL_{D_1}\ev^*=0=\mu^\de_*(i_\de\times\id)^*\ev^*,
\]
upon using \eqref{muLD}.
\end{proof}

{}From now on we assume $\omega=\dd\vartheta$.
As in Theorem~\ref{thm-Dinv}, we define the odd $1$\ndash form
$\bvartheta=\mu_*\ev^*\vartheta$. We also introduce the even $1$\ndash form
$\tau=\mu_*^\de\ev^*_\de\vartheta$.
Observe that $\tau$ is the obstruction to having $\Hat D$ Hamiltonian;
in fact, by Lemma~\ref{lem-ldhat} we have 
\begin{equation}\label{patronus}
\iota_{\Hat D}\bomega=
\iota_{\Hat D}\dd\bvartheta=
-\dd\iota_{\Hat D}\bvartheta+
\tau.
\end{equation}
In order to proceed we consider then the zero locus $\calC$
of $\tau$. This can also be described as the common zero locus
of a set of functions:
For every vector field $\xi$ on
$Y^{\Pi T\Sigma}$, we define $H_\xi:=\iota_\xi\tau$;
then $\calC=\{f\in Y^{\Pi T\Sigma}:H_\xi(f)=0\,\forall\xi\}$.
Observe that $H_\xi(f)$ actually depends only on the values that
$\xi(f)$, seen as an element of $\Gamma(\Pi T\Sigma, f^*TY)$, 
assumes on $\de\Sigma$.
If we denote by $Z(\vartheta)$ the zero locus of $\vartheta$,
we also obtain $\calC=\{f\in Y^{\Pi T\Sigma}:\iota_\de^*f(\Pi T\de\Sigma)
\subset Z(\vartheta)\}$.

The above construction turns out however to be too singular, for
no Hamiltonian vector fields correspond to
the functions $H_\xi$. We regularize then as follows.
Let $U$ be an open collar of $\de\Sigma$ in $\Sigma$.
We choose on $U$ a normal coordinate $x_{\rm n}$ and note by
$\theta_{\rm n}$ the corresponding odd counterpart to $x_{\rm n}$
on $\Pi TU$. Then we define
\[
\tau^U=\mu^U_*(\ev_U^*\vartheta\,\theta_{\rm n}),
\]
with obvious meaning of notations. 
Thus, for every vector field $\xi$ on $Y^{\Pi TU}$, we can write
\[
H_\xi^U(f):=\iota_\xi\tau^U=
\int_{\Pi TU}\xi^i(f,u)\vartheta_i(f(u))\theta_{\rm n}\;\mu(u).
\]
Let $\frX^U$ be the space of vector fields $\xi$ on $Y^{\Pi T\Sigma}$ 
with the property that $\xi(f)$, seen as an element of 
$\Gamma(\Pi T\Sigma, f^*TY)$, 
has support in $U$ for any $f\in Y^{\Pi T\Sigma}$. Any element of
$\frX^U$ determines a vector field on $Y^{\Pi TU}$ which, by abuse
of notation, we will denote by the same letter. So we can define
\begin{align*}
\calC^U&=\{f\in Y^{\Pi T\Sigma}:H_\xi^U(f)=0\,\forall\xi\in\frX^U\}\\
&=\{f\in Y^{\Pi T\Sigma}:\theta_{\rm n}f(\Pi TU)\subset Z(\vartheta)\}.
\end{align*}
We have then the following obvious
\begin{Lem}\label{lem-C}
For every two collars $U$ and $V$ with $U\subset V$,
one has $\calC^V\subset\calC^U$. Moreover, $\calC^U\subset\calC$
for every collar $U$.
\end{Lem}
Now observe that the functions $H_\xi^U$, $\xi\in\frX^U$,
define first-class constraints on $Y^{\Pi T\Sigma}$. In fact,
their Hamiltonian vector fields are proportional to $\theta_{\rm n}$.
So their application to any $H_\xi^U$ automatically vanishes
thanks to $(\theta_{\rm n})^2=0$. Thus, we can define a new 
$P$\ndash space by considering the quotient of $\calC^U$ by the foliation
generated by the Hamiltonian
vector fields of all $H_\xi^U$. We denote this reduced space by
$Y^{\Pi T\Sigma}//\tau^U$.
Then we have the following
\begin{Thm}\label{thm-dementor}
$\Hat D$ is Hamiltonian on $Y^{\Pi T\Sigma}//\tau^U$
with Hamiltonian function $\Hat S=-\iota_{\Hat D}\bvartheta$ satisfying
$\sBV{\Hat S}{\Hat S}=0$.

If moreover $Q$ is Hamiltonian and its Hamiltonian function
$S$ is locally constant on the zero locus of $\vartheta$, then
the restriction of $\Check S$ to $\calC^U$ is an invariant function.
So $\Check S$ descends to a function on $Y^{\Pi T\Sigma}//\tau^U$
that we will denote by the same symbol. Of course,
$\sBV{\Check S}{\Check S}=0$ still holds if $\BV SS=0$.

Finally, if we further assume that $S$ vanishes on the zero locus 
of $\vartheta$, then $\sBV{\Hat S}{\Check S}=0$.
In particular, $\sfS=\Hat S+\Check S$ satisfies the master equation
$\sBV\sfS\sfS=0$.
\end{Thm}
\begin{proof}
The first statement follows from \eqref{patronus} and Lemma~\ref{lem-C}.
We have first to check that $\Hat S$ restricted to $\calC^U$
is invariant. This is a consequence of the identities
\[
\sBV{\Hat S}{H_\xi^U}=\LL_{\Hat D}\iota_\xi\tau^U=
\LL_{\Lie{\Hat D}\xi}\tau^U\pm\iota_\xi\LL_{\Hat D}\tau^U=0.
\]
The last equality holds on $\calC^U$
since $\Lie{\Hat D}\xi\in\frX^U$ for
$\xi\in\frX^U$ and because the restriction of 
$\ev^*_U\vartheta\,\theta_n$ to $\Pi T\de U$ vanishes (as one has to
set $\theta_{\rm n}=0$).
Finally, $\sBV{\Hat S}{\Hat S}=-\LL_{\Hat D}\iota_{\Hat D}\bvartheta
=\iota_{\Hat D}\LL_{\Hat D}\bvartheta=\iota_{\Hat D}\tau=0$.

As for the second statement,
let us denote by $\theta_{\rm n} X$ the Hamiltonian
vector field of $H_\xi^U$. Observe that $X(f,\bullet)$ has support
in $U$, so
\[
\sBV{H_\xi^U}{\Check S}(f)=
\int_{\Pi TU}\theta_{\rm n}X^i(f,u)\de_i S(f(u))\;\mu.
\]
However, if we restrict to $\calC^U$, then $\theta_{\rm n}f(u)$ 
belongs to the zero locus of $\vartheta$ for all $u\in U$, and
$S$ is constant on it by assumption.
So $\sBV{H_\xi^U}{\Check S}=0$.

Since the vector
fields $\Hat D$ and $\Check Q$ commute, it follows immediately
that the bracket $\sBV{\Hat S}{\Check S}$ 
is a constant function, but we want this
function to vanish identically.
Actually,
by Lemma~\ref{lem-ldhat} we have
\[
\sBV{\Hat S}{\Check S}=
\LL_{\Hat D}\Check S=\mu^\de_*\ev_\de^*S=0,
\]
the last equality holding on $\calC$ (and so on $\calC^U$) since
$f(u)$ must belong to the zero locus of $\vartheta$ for $u\in\de\Sigma$,
and $S$ vanishes on it by assumption.
\end{proof}

In the application of the above construction to the perturbative
evaluation of a path integral, the idea is to define it first on
$Y^{\Pi T\Sigma}//\tau^U$ for a given collar $U$.
Next one should compute the propagators
and the vertices and finally consider the limit for
$U$ shrinking to $\de\Sigma$. 

Another possibility, see e.g.\ next
section, is to fix representatives in $Y^{\Pi T\Sigma}//\tau^U$
as elements of $\calC^U$ satisfying some extra conditions.
Let us denote by $\Tilde\calC^U$ this further constrained space
isomorphic to $Y^{\Pi T\Sigma}//\tau^U$. If $U\subset V$ implies
$\Tilde\calC^V\subset\Tilde\calC^U$, we define $\Tilde\calC$
as $\bigcup\Tilde\calC^U$. Then the path integral may be defined
on a Lagrangian submanifold (see Remark~\ref{rem-BVquant} in the
next subsection) of $\Tilde\calC$.

\subsubsection{Remarks}
\begin{Rem}[dependency on $\vartheta$]
In the general case, described in Theorem~\ref{thm-Dinv},
the choice of a potential $\vartheta$ for the symplectic form
$\omega$ is irrelevant. In fact, $\vartheta$ enters only in the
definition of $\Hat S$. But if we choose $\vartheta'=\vartheta+\dd\sigma$,
we obtain 
$\Hat S':=-\iota_{\Hat D}\bvartheta'=\Hat S-\LL_{\Hat D}\mu_*\ev^*\sigma$, 
and this equals $\Hat S$ by Lemma~\ref{lem-Dinv}.

In the case of $X=\Pi T\Sigma$, $\de\Sigma\not=\emptyset$, the choice of
$\vartheta$ is instead
essential as it determines the boundary conditions
through $\tau$. As different but cohomologous $\vartheta$'s may have
different zero sets, this also affects which functions $S$ on $Y$
are allowed by Theorem~\ref{thm-dementor}. (Think, e.g., of 
$\Pi T^*\bbR^m$ with $\omega=\dd p_i\dd q^i$ and
$\vartheta=a p_i\dd q^i+(a-1) q^i\dd p_i$ for different choices of
$a$.)
\end{Rem}
\begin{Rem}[ghost number]\label{rem-gn}
As observed in~\ref{ssec-zgrad}, in quantum field theory one is
interested in the full $\bbZ$\ndash grading. So, in order to have
$\bomega$ of ghost number $-1$, we must choose $\omega$ to have 
ghost number $n-1$, if $X$ is of type $(m,n)$. Moreover, we must
assume that $Q$ has ghost number one, which implies that its
Hamiltonian function $S$ must have ghost number $n$.
In the case $X=\Pi T\Sigma$, we require then $\gh\omega=s-1$ and
$\gh S=s$ with $s=\dim\Sigma$.
\end{Rem}
\begin{Rem}[twisted supersymmetry]\label{rem-SUSY} 
Since the canonical measure $\mu$ on $\Pi T\Sigma$ is $K_v$\ndash invariant
for every vector field $v$ on $\Sigma$, we conclude by Lemma~\ref{lem-Dinv}
that $\bomega$, $\bvartheta$ and $\Check S$ are
$\Hat K_v$\ndash invariant. (Actually, $\Hat K_v$ is even Hamiltonian
with Hamiltonian 
function $-\iota_{\Hat K_v}\bvartheta$.)
We obtain then
\[
\Lie{\Hat K_v}\sfQ=\Lie{\Hat K_v}{\Hat D}=\Hat L_v.
\]
In particular, we may locally choose constant vector fields 
$v_\mu:=\frac{\de}{\de u^\mu}$. 
We then have $L_{v_\mu}=\frac{\de}{\de u^\mu}$. This means that
$\sfQ$ and $\Hat K_{v_\mu}$, $\mu=1,\dots,\dim\Sigma$, are generators
of a twisted supersymmetry.
\end{Rem}
\begin{Rem}[classical theory]\label{rem-cl}
The AKSZ method induces {\em a posteriori}\/ a classical theory of which
the above discussion is the BV version.
Distinguish the components of a map $f\colon X\to Y$ according to
the ghost number. The components of nonnegative ghost number are called
the {\sf fields}, the others the {\sf antifields}. Among the fields one further
distinguish between the {\sf classical fields} (zero ghost number) and the
{\sf ghosts} (positive ghost numbers). The {\sf classical action}
$S^{\rm cl}$ is
obtained by setting all the antifields in $\sfS$ to zero; it depends
only on the classical fields. The action of the $Q$-vector field
on the fields at zero antifields is called the {\sf BRST operator} 
and generates
the infinitesimal symmetries of the theory.
\end{Rem}
\begin{Rem}[BV quantization]\label{rem-BVquant}
Given a solution of the master equation,\footnote{Actually, 
one needs a solution of the so-called 
quantum master equation
\[
\sBV\sfS\sfS-2\ii\hbar\Delta\sfS=0,
\]
where $\Delta$ is an operator that depends on the measure in the functional
integral and so is well-defined only after a regularization has been 
chosen.}
the rules of the game for defining (in perturbation theory) the
functional integral of $\exp\ii S^{cl}/\hbar$ over the classical fields
are as follows: Choose a Lagrangian submanifold ({\sf gauge fixing})
of the space of fields and antifields, and integrate over it
the new weight $\exp\ii \sfS/\hbar$. This procedure works whenever
the new action is nondegenerate on the chosen Lagrangian submanifold.
Actually, the AKSZ construction has not enough room for this step, but one
just has to enlarge the space $Y^X$ to include enough
Lagrange multipliers ($\lambda$)
and antighosts ($\bar c$)---together with their canonically conjugated 
variables $\lambda^\dagger$ and $\bar c^\dagger$---in order to 
gauge-fix all symmetries in the extended action given by
$\sfS$ plus terms of the form $\bar c^\dagger\lambda$.
\end{Rem}

\begin{Rem}[Applications]
The AKSZ method was foreshadowed by Witten \cite{W2} in the sigma-model
interpretation of Chern--Simons theory \cite{W1}. 
In \cite{AKSZ} it was explained that Chern-Simons theory enters this scheme
by choosing $\Sigma$ a $3$\ndash manifold and $Y=\frg[1]=\Pi\frg$, with $\frg$
a metric Lie algebra; in the same paper 
the method was also applied to the A and B 
models \cite{W3}. In \cite{K-1} the application to the Rozansky--Witten 
\cite{RW} theory was considered.
In \cite{CR} the method has been applied to $BF$ theories
by choosing $\Sigma$ an $s$\ndash manifold and $Y=\frg[1]\times\frg^*[s-2]$.
In \cite{P}, topological open branes are defined by choosing
$\Sigma$ an $s$\ndash manifold with boundary and $Y=(\Pi T)^{s-2}\Pi T^*M$,
with $M$ an ordinary manifold.
\end{Rem}

\section{The Poisson sigma model}\label{s-3}
Let $\Sigma$ and $M$ be ordinary manifolds. We want to 
apply the AKSZ construction described in the previous Section by choosing
$X=\Pi T\Sigma$ and $Y=\Pi T^* M$, with the $Q$\ndash vector field
$D$, the canonical measure $\mu$ and the odd symplectic structure
$\omega$ introduced in Examples~\ref{exa-ptsigma}, 
\ref{exa-mu} and \ref{exa-ptm} respectively. Since $\omega$ has
in this case ghost number $1$ we will consider $\Sigma$ to be
$2$\ndash dimensional, according to the discussion in
Remark~\ref{rem-gn}.

\subsection{The $P$-structure}
An element of $Y^X$ in the present case is a pair $(\sfX,\sfeta)$ where
$\sfX$ is a map $\Pi T\Sigma\to M$ and $\sfeta$ is a section of
$\sfX^*\Pi T^*M$. In local coordinates on 
$\Pi T\Sigma$ (see Example~\ref{exa-ptsigma}), we may write
\begin{align*}
\sfX &= X + \theta^\mu\eta^+_\mu-\frac12\theta^\mu\theta^\nu\beta^+_{\mu\nu},\\
\sfeta &= \beta + \theta^\mu\eta_\mu+\frac12\theta^\mu\theta^\nu X^+_{\mu\nu},
\end{align*}
with $X$ an ordinary map $\Sigma\to M$ and
\begin{align*}
\eta^+ &\in\Gamma(\Sigma,T^*\Sigma\otimes X^*TM)\otimes\Lambda_-[-1],\\
\beta^+ &\in\Gamma(\Sigma,\wedge^2T^*\Sigma\otimes X^*TM)
\otimes\Lambda_+[-2],\\
\beta &\in \Gamma(\Sigma,X^*T^*M)\otimes\Lambda_-[1],\\
\eta &\in \Gamma(\Sigma,T^*\Sigma\otimes X^*T^*M)\otimes\Lambda_+[0],\\
X^+ &\in \Gamma(\Sigma,\wedge^2T^*\Sigma\otimes X^*T^*M)\otimes\Lambda_-[-1],
\end{align*}
where $\Lambda_\pm$ are the odd and even parts of $\Lambda$, and
we follow the notations of~\ref{ssec-zgrad}.
The classical fields (see Remark~\ref{rem-cl}) are then $X$ and $\eta$,
the only ghost is $\beta$, and the other components are antifields. 
By the symplectic form $\bomega$, the superfields $\sfX$ and $\sfeta$ 
are canonically conjugated. In terms of their components, each field
is canonically conjugated to the corresponding antifield with an upper $+$.

The canonical symplectic form is exact: $\omega=\dd\vartheta$. We choose
$\vartheta=p_i\dd x^i=\braket p{\dd x}$,
where $\braket{\ }{\ }$ is the canonical pairing of vectors and covectors 
on $M$, and define $\bvartheta=\mu*\ev^*\vartheta$ as in the previous
section. To give a more explicit description of $\bvartheta$, we observe
that the tangent bundle of $\Pi T^*M^{\Pi T\Sigma}$ has a splitting
$T(\Pi T^*M^{\Pi T\Sigma})=\frX\oplus\frE$ given by vectors
``in the directions of $\sfX$ and in the directions of $\sfeta$''
respectively; viz.,
\begin{align*}
\frX_{(\sfX,\sfeta)} 
&= T_{(\sfX,\sfeta)}M^{\Pi T\Sigma}\simeq\Gamma(\Pi T\Sigma,X^*TM),\\
\frE_{(\sfX,\sfeta)} 
&= T_{(\sfX,\sfeta)}\Gamma(\Pi T\Sigma,X^*\Pi T^*M)\simeq
\Gamma(\Pi T\Sigma,\sfX^*\Pi T^*M).
\end{align*}
Then we can write
\[
\bvartheta(\sfX,\sfeta)(\bxi\oplus\sfe)=\int_{\Pi T\Sigma} 
\braket\bxi\sfeta\;\mu,
\]
for any
$\bxi\in\frX_{(\sfX,\sfeta)}$ and $\sfe\in\frE_{(\sfX,\sfeta)}$.

\subsection{The $QP$-structure and the classical action}
As in the previous section we define a $Q$\ndash structure on $Y^X$
starting from $Q$\ndash structures on $X$ and on $Y$. 
On $Y=\Pi T^* M$ we look for a Hamiltonian vector field of ghost
number $1$. Since $\omega$ has ghost number $1$, this means that we have
to look for a Hamiltonian function $S$ of ghost number $2$, i.e., for
a function corresponding to a bivector 
field\footnote{\label{f-multi}More generally, 
we may consider any multivector field such that
the component of the term of order $k$ has ghost number $k-2$.}
as in Example~\ref{exa-ptmPQ};
viz., we define $S_\alpha$ and $Q_\alpha$ 
as in \eqref{SQ}. We recall that a manifold
endowed with a bivector field $\alpha$ satisfying
$\Lie\alpha\alpha=0$ is called a Poisson manifold, for $\alpha$
defines a Poisson bracket on the algebra of functions.
By \eqref{checks}, it follows that
the Hamiltonian function of the vector field $\Check Q_\alpha$ is
\[
\Check S_\alpha=\frac12\int_{\Pi T\Sigma}(\alpha\circ\sfX)(\sfeta,\sfeta)\;\mu.
\]
On $X=\Pi T\Sigma$
we take the $Q$\ndash vector field $D$ of Example~\ref{exa-ptsigma}, which
has ghost number $1$. As in Theorem~\ref{thm-dementor}, 
we may define $\Hat S=-\iota_{\Hat D}\bvartheta$ getting
\[
\Hat S= -\int_{\Pi T\Sigma}\braket{D\sfX}\sfeta\;\mu,
\]
where $D\sfX$ is a short-hand notation for the projection of $\Hat D$
to $\frX$.

By Theorem~\ref{thm-Dinv}, the function $\sfS_\alpha=\Hat S+\Check S_\alpha$
satisfies the master equation if $\de\Sigma=\emptyset$. As described
in Remark~\ref{rem-cl}, we can recover the classical theory to which
$\sfS$ corresponds by setting the antifields to zero; we get then the
classical action
\[
S^{\rm cl}_\alpha=\int_\Sigma \left(
\braket\eta{\dd X}+\frac12(\alpha\circ X)(\eta,\eta)
\right),
\]
which describes the so-called Poisson sigma model \cite{Ik, SS}.

We now turn to discuss the case when $\Sigma$ has a boundary. 
Observe that the zero locus of $\vartheta$
is in the present case the zero section of $\Pi T^*M$, i.e., $p^i=0$.
So our $S$ satisfies\footnote{If we work in the more general setting of
footnote~\ref{f-multi}, we should however exclude from $\alpha$
any term of order zero, i.e., functions. 
}
the strongest assumption of Theorem~\ref{thm-dementor} and
we get a solution of the master equation for every collar $U$.
We next describe $\calC^U$. 
First observe that the even $1$\ndash form
$\tau^U$ vanishes when contracted with
vectors in $\frE$. 
If instead $\bxi$ is
a vector in $\frX$---i.e., 
$\bxi(\sfX,\sfeta)\in\Gamma(\Pi T\Sigma,X^*TM)$---, then we have
\[
H_\bxi^U=\iota_\bxi\tau^U=
\int_{\Pi TU} \braket\bxi\sfeta\theta_{\rm n}\;\mu^U.
\]
Thus, $\calC^U$ is the space of maps $(\sfX,\sfeta)$ satisfying
$\theta_{\rm n}\iota_U^*\sfeta=0$. The Hamiltonian vector field
$H_\bxi$ corresponding to a constant vector field $\bxi$
is also easily computed to be  $\bxi\theta_{\rm n}$.
This means that, on $\Pi T^* M^{\Pi T\Sigma}//\tau^U$,
$\sfX$ is defined up to translations in the normal direction
supported in $U$.
We can fix a representative by choosing any direction ${\rm n}'$
transversal to $\de\Sigma$ (e.g., ${\rm n}'={\rm n}$) and
by requiring that the ${\rm n}'$\ndash components of $\sfX$ vanish.
This means $\beta^+=0$ and $\eta^+_{\rm n'}=0$.
We can finally remove the regularization by considering
\begin{equation}\label{bc}
\Tilde\calC=\{(\sfX,\sfeta):\beta(u)=0,\ \eta_{\rm t}(u)=0,\ 
\eta^+_{\rm n'}(u)=0,\ \beta^+_{\rm n't}(u)=0\ \forall u\in\de\Sigma\},
\end{equation}
where $\rm t$ denotes the direction tangent to $\de\Sigma$.

The above 
discussion leads exactly to the BV action of \cite{CF} with the same
boundary conditions. We refer to it for the introduction of antighosts
and Lagrange multipliers, for the discussion of the quantum master
equation and for the relation between the Poisson sigma model on
the disk and Kontsevich's formula \cite{K} for the deformation quantization
 of $M$.

\subsection{Target diffeomorphisms}
We conclude with a brief discussion on the effect of a diffeomorphism $\phi$
of $M$ on the whole theory. We denote by $\Phi$ the canonical
extension of $\phi$ to $\Pi T^*M$:
\[
\Phi(x,p)=(\phi(x),(\dd\phi(x)^*)^{-1}p).
\]
It leaves invariant the canonical symplectic form $\omega$ as well
as the canonical $1$\ndash form $\vartheta$. Moreover, one has
$\Phi^*S_{\phi_*\alpha}=S_\alpha$.
Let $\Check\Phi$ be the 
corresponding right action on $\Pi T^*M^{\Pi T\Sigma}$:
\[
\Check\Phi(f)=\Phi\circ f, \qquad f\in\Pi T^*M^{\Pi T\Sigma}.
\]
It follows that $\Check\Phi$ is a symplectomorphism of 
$(\Pi T^*M^{\Pi T\Sigma},\bomega)$ with the property
$\Check\Phi^*\Check S_{\phi_*\alpha}=\Check S_\alpha$.
Moreover, since $\Hat D$ corresponds
to an (infinitesimal) left action and $\Check\Phi^*\bvartheta=\bvartheta$, 
we also obtain that $\Hat S$
is $\Check\Phi$\ndash invariant. So we conclude that
\begin{equation}\label{PhiS}
\Check\Phi^*\sfS_{\phi_*\alpha}=\sfS_\alpha.
\end{equation}
This means that the Poisson sigma model
is invariant under target diffeomorphisms at the classical level, and it 
is so at the quantum level if the measure on the space of fields
is invariant.

To show that this is so at the quantum level, we restrict for simplicity 
to infinitesimal diffeomorphisms.
Let $\xi$ be a vector field on $M$ and $\Xi$ the corresponding vector
field on $\Pi T^*M$. It is not difficult to see that $\Xi$ has
Hamiltonian $S_\xi(x,p)=\xi^i(x) p_i$. Let $\Check\Xi$ be the
corresponding vector field on $\Pi T^*M^{\Pi T\Sigma}$.
It is a Hamiltonian vector field with Hamiltonian 
$\Check S_\xi=\mu_*\ev^*S_\xi$. 
We recall that a canonical transformation in the BV formalism extends to 
the quantum level
(i.e., it preserves the path integral measure), if it is in the kernel
of the BV Laplacian (which is constructed using the measure).
We proved in \cite{CF} that this is the case for a function like
$\Check S_\xi$ (actually, with $\xi$ any multivector field)
in reasonable regularizations.

It is however well-known that Kontsevich's formula for the $\star$-product
does not transform well under diffeomorphisms. We would like to end the 
paper with some comment on this apparent contradiction.
First observe that the infinitesimal form of \eqref{PhiS},
\[
\LL_{\Check\Xi}\sfS_\alpha=\sfS_{\Lie\xi\alpha},
\]
may be obtained by the identities
\begin{gather*}
\LL_{\Check\Xi}\Check S_\alpha = \sBV{\Check S_\xi}{\Check S_\alpha}=
\Check S_{\Lie\xi\alpha},\\
\LL_{\Check\Xi}\Hat S = -\sBV{\Hat S}{\Check S_\xi}=
-\Hat D\Check S_\xi=-\mu^\de_*\ev_\de^*S_\xi=0.
\end{gather*}
Now the last equality follows from the fact that $\ev_\de^*S_\xi$
vanishes upon using the boundary conditions in \eqref{bc}. However,
in the definition of the $\star$\ndash product from the path integral
of the Poisson sigma model, one has to insert boundary observables
(i.e., evaluate some functions on $M$ at points $X(u)$ with 
$u\in\de\Sigma$), see \cite{CF}. These introduce singularities that have
to be removed by choosing a small neighborhood $U$ of the boundary
points of insertions and
integrating on $\Pi T(\Sigma\setminus U)$ instead of $\Pi T\Sigma$.
As a consequence, in the last equation the term 
$S^C_\xi=\int_{\Pi T C} \ev^*S_\xi\;\mu^C$ survives, where $C$ is the boundary
of $U$ modulo the boundary of $\Sigma$. This spoils the original invariance
of the Poisson sigma model under target diffeomorphisms.
It also suggests that the deformed action of a diffeomorphism
on a function should be given by the expectation value of $S^C_\xi$
(i.e., $\sum_n(\ii\hbar)^n \calU_{n+1}(\xi,\alpha,\dots,\alpha)/n!$,
where $\calU$ denotes Kontsevich's $L^\infty$\ndash morphism \cite{K}).
This idea, also related to observations in \cite{L}, 
has been used in \cite{CFT}
to clarify the globalization of Kontsevich's formula.

\end{document}